\documentclass[sigconf,screen,final,nonacm]{acmart}
\AtBeginDocument{%
  }

\usepackage{algorithm}
\usepackage{algpseudocode}
\algblockdefx[PARFOR]{PARFOR}{ENDPARFOR}%
  [1]{\textbf{Parallel for} #1 \textbf{do}}%
      {\textbf{end for}}

\begin{document}

\title{A High Performance GPU CountSketch Implementation and Its Application to Multisketching and Least Squares Problems}

\author{Andrew J. Higgins}
\email{ajhiggi@sandia.gov}
\affiliation{%
  \institution{Sandia National Laboratories}
  \city{Albuquerque}
  \state{New Mexico}
  \country{USA}
}

\author{Erik G. Boman}
\email{egboman@sandia.gov}
\affiliation{%
  \institution{Sandia National Laboratories}
  \city{Albuquerque}
  \state{New Mexico}
  \country{USA}
}

\author{Ichitaro Yamazaki}
\email{iyamaza@sandia.gov}
\affiliation{%
  \institution{Sandia National Laboratories}
  \city{Albuquerque}
  \state{New Mexico}
  \country{USA}
}

\renewcommand{\shortauthors}{Higgins et al.}

\begin{abstract}
Random sketching is a dimensionality reduction technique that approximately preserves norms and singular values up to some $O(1)$ distortion factor with high probability. The most popular sketches in literature are the Gaussian sketch and the subsampled randomized Hadamard transform, while the CountSketch has lower complexity. Combining two sketches, known as multisketching, offers an inexpensive means of quickly reducing the dimension of a matrix by combining a CountSketch and Gaussian sketch.

However, there has been little investigation into high performance CountSketch implementations. In this work, we develop an efficient GPU implementation of the CountSketch, and demonstrate the performance of multisketching using this technique. We also demonstrate the potential for using this implementation within a multisketched least squares solver that is up to $77\%$ faster than the normal equations with significantly better numerical stability, at the cost of an $O(1)$ multiplicative factor introduced into the relative residual norm.
\end{abstract}

\begin{CCSXML}
<ccs2012>
   <concept>
       <concept_id>10002950.10003648.10003671</concept_id>
       <concept_desc>Mathematics of computing~Probabilistic algorithms</concept_desc>
       <concept_significance>300</concept_significance>
       </concept>
   <concept>
       <concept_id>10002950.10003705.10003707</concept_id>
       <concept_desc>Mathematics of computing~Solvers</concept_desc>
       <concept_significance>300</concept_significance>
       </concept>
   <concept>
       <concept_id>10002950.10003705.10011686</concept_id>
       <concept_desc>Mathematics of computing~Mathematical software performance</concept_desc>
       <concept_significance>300</concept_significance>
       </concept>
 </ccs2012>
\end{CCSXML}

\ccsdesc[300]{Mathematics of computing~Probabilistic algorithms}
\ccsdesc[300]{Mathematics of computing~Solvers}
\ccsdesc[300]{Mathematics of computing~Mathematical software performance}

\keywords{Randomized Linear Algebra, GPU Performance}

\maketitle

\section{Introduction to Sketching} \label{sec:intro}
Random sketching--the process of applying a random linear operator to a matrix/vector in order to reduce its dimension while approximately preserving its key information--has become ubiquitous throughout randomized numerical linear algebra. The most popular type of sketches are subspace embeddings \cite{martinsson_tropp_2020}, which are specific sketches that preserve inner products among all vectors in a fixed dimensional subspace with high probability. 

\begin{definition}[$\varepsilon$-subspace embedding] \label{def:subspaceEmbedding}
Given $\varepsilon \in [0,1)$, the sketch operator $S \in \mathbb{R}^{d} \rightarrow \mathbb{R}^k$ is an \emph{$\varepsilon$-subspace embedding} for the subspace $\mathcal{V} \subset \mathbb{R}^d$ if $\forall x, y \in \mathcal{V}$, 
\begin{equation}
    | \langle x,y \rangle - \langle Sx,Sy \rangle | \leq \varepsilon \|x\|_2\|y\|_2. \label{eq:innerProdEmbedding}
\end{equation}
\end{definition}

\begin{definition}[$(\varepsilon, \delta, n)$ oblivious $\ell_2$-subspace embedding] \label{def:obliviousSubspaceEmbedding}
$S \in \mathbb{R}^{d} \rightarrow \mathbb{R}^k$ is an \emph{$(\varepsilon, \delta, n)$ oblivious $\ell_2$-subspace embedding} if it is an $\varepsilon$-subspace embedding for any fixed $n$-dimensional subspace $\mathcal{V} \subset \mathbb{R}^d$ with probability at least $1-\delta$.
\end{definition}

The most famous and popular type of subspace embedding used in the literature is the Gaussian sketch, where $S \in \mathbb{R}^{k \times d}$ and $s_{i,j} \sim \mathcal{N}(0,k^{-1})$ for $k = O\left(\frac{n - \log(\delta)}{\epsilon^2}\right)$ \cite{Balabanov:2019}. Specifically, $k = n/\varepsilon^{-2}$ ensures $S$ forms a $\varepsilon$-subspace embedding with high probability \cite{martinsson_tropp_2020}. However, since its action on a vector/matrix involves multiplication with a dense matrix, its application to $A \in \mathbb{R}^{d \times n}$ has a computational complexity of $O(dn^2)$. 

The cheapest known sketch is the CountSketch \cite{CountSketchOriginal}, which can be interpreted as $S \in \mathbb{R}^{k \times d}$, where $S$ has exactly one $\pm 1$ per column, chosen uniformly at random \cite{WoodruffSketching}. While this can be applied to $A \in \mathbb{R}^{d \times n}$ in only $O(dn)$ operations, the drawback of the CountSketch is that one needs to choose $k = O\left( \frac{n^2}{\varepsilon^2\delta} \right)$ to guarantee that it will be a $(\varepsilon, \delta, n)$ subspace embedding \cite{CountSketchSize}. 

Another subspace embedding that is popular in literature for its theoretical properties is the Subsampled Randomized Hadamard Transform (SRHT) \cite{Tropp11}. This method is based on applying a sign flip, a fast Walsh-Hadamard transform (FWHT) \cite{FWHT}, and a row sampling operation. The SRHT offers a compromise between the Gaussian sketch and the CountSketch in terms of both the embedding dimension and computational complexity, requiring an embedding dimension $k = O(n \log n)$ to ensure the subspace embedding properties hold and a complexity as low as $O(dn \log k)$ operations \cite{Ailon09, Tropp11}. However, in practice, the SRHT preserves the subspace embedding properties provided $k = O(n)$ \cite{martinsson_tropp_2020}. The drawback of the SRHT is that it requires a good implementation of the FWHT, which is non-trivial on a GPU and needs to access $A \in \mathbb{R}^{d \times n}$ $O(\log k)$ times.

Instead of using the SRHT, which may be bottlenecked by the FWHT implementation, one can take advantage of the CountSketch's low cost while still reducing the dimension further than $O(n^2)$ by using a multisketching strategy. In its most generic form, one chooses two sketches $S_1 \in \mathbb{R}^{k_1 \times d}$ and $S_2 \in \mathbb{R}^{k_2 \times k_1}$ to a vector $x \in \mathbb{R}^d$ via $S_2(S_1 x)$. Specifically, $S_1$ is chosen so that it can quickly reduce the dimension of $x$, such as via CountSketch, and then $S_2$ is chosen to reduce the dimension further to a more manageable size. The most straigthforward example of this approach is the so-called ``Count-Gauss" sketch \cite{Kapralov_Potluru_PMLR_16}, whose efficiency has implicitly been demonstrated on a GPU while computing a randomized QR factorization, where $S_1x$ is computed using a sparse matrix-vector multiply  \cite{higgins2024}. We list the properties of each of the aforementioned sketching techniques in Table \ref{tab:sketchProperties}.

\begin{table}[h]
\centering
\scriptsize
\begin{tabular}{l||l|l|l|l}
                      & Embed Dim. & Arithmetic & Read/Writes & Max Distortion \\
\hline                
Gaussian & $\epsilon^{-2}n$ & $dn^2$ & $dn$ & $1 + \epsilon$\\
SRHT & $\epsilon^{-2}n \log n$ & $dn \log n$ & $d n \log n$ & $1 + \epsilon$ \\
CountSketch & $\epsilon^{-2}n^2$ & $dn$ & $dn$ & $1 + \epsilon$ \\
MultiSketch$(\epsilon_1, \epsilon_2)$ & $\epsilon_2^{-2}n$ & $dn + n^4$ & $dn + n^4$ & $(1 + \epsilon_1)(1+\epsilon_2)$ 
\end{tabular}
\caption{Asymptotically optimal embedding dimensions, arithmetic, and read/writes required to apply each of the sketching methods to a dense matrix $A \in \mathbb{R}^{d \times n}$. `Multisketch$(\epsilon_1, \epsilon_2)$' denotes a CountSketch with distortion $\epsilon_1$ followed by a Gaussian with distortion $\epsilon_2$.}
\label{tab:sketchProperties}
\end{table}

\section{Sketching for Least Squares Problems}

One of the earliest and most popular applications of random sketching was its integration into least squares problems, using a ``sketch-and-solve" approach \cite{Drineas2006, Sarlos:2006}. Specifically, suppose we have an overdetermined least squares problem $\min_{x} \| Ax - b \|_2$ where $A \in \mathbb{R}^{d \times n}$, $b \in \mathbb{R}^d$, $d \gg n$, and let $S \in \mathbb{R}^{k \times d}$ be a random sketch operator. The ``sketch-and-solve" approach approximately solves the problem by sketching both the coefficient matrix $A$ and the right hand side vector $b$ and solving the resulting reduced least squares problem $\min_{x} \| Sb - SAx \|_2$ instead, and is written explicitly in algorithm \ref{alg:sketchNSolve}. Provided the sketch operator $S$ is an $(\epsilon, \delta, n)$ oblivious $\ell_2$-subspace embedding, $\sqrt{1-\varepsilon} \| b - Ax \|_2 \leq \| S(b-Ax) \|_2$, $\forall x \in \mathbb{R}^n$. Thus, if $x_s = \text{argmin} \| S(b-Ax) \|_2$ and $x_t = \text{argmin} \| b - Ax \|_2$, 
\begin{align*}
    \|b - Ax_t \|_2 &\leq \|b - Ax_s \|_2 \leq \frac{1}{\sqrt{1-\varepsilon}} \| Sb - SAx_s \|_2 \\ &\leq \frac{1}{\sqrt{1-\varepsilon}} \| Sb - SAx_t \|_2 \leq \sqrt{\frac{1+\varepsilon}{1-\varepsilon}} \| b - Ax_t \|_2.
\end{align*}
Hence, the ``sketch-and-solve" solution $x_s$ gives an approximate least squares solution to the original problem $\min \|b - Ax\|_2$, up to some $O(1)$ distortion factor in $\left[1, \sqrt{\frac{1+\varepsilon}{1-\varepsilon}} \right]$. The distinct advantage of the ``sketch-and-solve" approach is that if one is able to cheaply apply the sketch operator and the sketch operator sufficiently reduces the dimension of the coefficient matrix, the dominant costs (sketching and matrix factorization) are cheap and one can approximately solve the least squares problem very quickly.

\begin{algorithm}[h]
         \begin{algorithmic}[1]  
             \vspace{0cm}
             \Statex \textbf{Input:}\phantom{aa} Matrix $A \in \mathbb{R}^{d \times n}$, vector $b \in \mathbb{R}^d$, sketch operator \phantom{aaaaaaaaa}$S: \mathbb{R}^d \rightarrow \mathbb{R}^k$
            \Statex {\textbf{Output:} Solution vector $x \in \mathbb{R}^n$}
             \vspace{0cm}
             \State Sketch $Y = SA$, $z = Sb$
             \State Compute economy QR factorization $[Q,R] = \text{qr}(Y,0)$
             \State Solve $x = R^{-1}\left(Q^Tz\right)$
         \end{algorithmic}
         \caption{Sketch-and-Solve Least Squares} \label{alg:sketchNSolve}
 \end{algorithm}

\section{Contributions}

In spite of the potential of fast sketching techniques to allow the sketch-and-solve approach to outperform the normal equations in theory, little attention has been dedicated to actually achieving this on modern hardwares, such as GPUs. The sketching algorithms most likely to allow the sketch-and-solve approach to outperform the normal equations, such as a multisketch approach, depend on a high performance CountSketch implementation, since it can be used to reduce the dimension of the original problem with the lowest complexity and theoretically only requires a single pass over the coefficient matrix and right hand side vector. Most CountSketches investigated in the randomized linear algebra literature use a simple sparse matrix multiply (SpMM or SpMV). However, SpMV/SpMM using a matrix with a random sparsity structure like the CountSketch on a vendor-provided library will likely lead to underwhelming performance, failing to maximize its potential. 

In this paper, we develop a straightforward, yet much higher performing GPU implementation of the CountSketch in CUDA and demonstrate its performance compared to using cuSPARSE to apply it. In Section \ref{sec:SketchPerf}, we demonstrate that using multisketching with this CountSketch implementation can be faster than computing the Gram matrix of $A \in \mathbb{R}^{d \times n}$. Hence, in Section \ref{sec:lsqr}, we show that a multisketched ``sketch-and-solve" approximate least squares solver is often faster than computing the solution using the normal equations, while improving stability significantly at the cost of an $O(1)$ distortion factor in the least squares residual. Finally, in Section \ref{sec:lsqr} we also compare the accuracy of the least squares solution achieved using a sketch-and-solve approach with the CountSketch compared to a rand\_cholQR based solver. Although the CountSketch can introduce a substantial distortion factor in theory, we observe that the distortion introduced by the sketch-and-solve approach with the CountSketch is negligible while being much faster than the rand\_cholQR based solver and only slightly slower than solving the normal equations.

\section{GPU CountSketch Implementation} \label{sec:countsketch}

We begin by revisiting an algebraically simple explicit definition of the CountSketch.
\begin{definition}[CountSketch] \label{def:countSketch}
    The $j^{th}$ column of the CountSketch $S \in \mathbb{R}^{k \times d}$ can be defined column-wise via $s_j = \sigma_j e_{r_j}$, where $\{\sigma_{\ell}\}_{\ell = 1}^d$ are i.i.d.~Rademacher random variables, and $e_m$ is the $m^{th}$ standard basis vector of $\mathbb{R}^k$, and $\{r_\ell\}_{\ell = 1}^d$ are i.i.d.~integer random variables drawn from the uniform distribution over $\{1, \dots, k\}$.
\end{definition}

Let $S \in \mathbb{R}^{k \times d}$ a CountSketch and $A \in \mathbb{R}^{d \times n}$. Clearly $Y = SA$ can be computed using a sparse-matrix times dense-matrix product. Alternatively, observe that for each $m \in \{1, \dots, k\}$
\begin{equation}
    Y_{m,:} = \sum_{\substack{1 \leq j \leq d \\ r_j = m}} \sigma_j A_{j,:},
\end{equation}
where $A_{j,:}$ denotes the $j^{th}$ row of $A$. Observe that multiplying by the random signs $\sigma_j$ and adding is equivalent to either adding or subtracting $A_{j,:}$, depending on its sign. Hence, the scalar multiply by $\sigma_j$ can be avoided and just be controlled by a boolean variable $s_j$ instead. Additionally, if one wants to execute this in a parallel reduction, since several rows of $A$ may simultaneously be added to a single row of $Y$, one should add rows atomically. Provided the random signs $\{s_{\ell}\}_{\ell = 1}^d$ and row mappings $\{t_\ell\}_{\ell = 1}^d$ are known, there is a simple parallel algorithm for implementing the CountSketch.

\begin{algorithm}[h]
         \begin{algorithmic}[1]  
             \vspace{0cm}
             \Statex \textbf{Input:}\phantom{aa} Matrix $A \in \mathbb{R}^{d \times n}$, integer vector $r \in \{1, \dots, k\}^d$, \phantom{aaaaaaaaaa}length $d$ boolean vector $s$ 
            \Statex {\textbf{Output:} Sketched matrix $Y \in \mathbb{R}^{k \times d}$}
             \vspace{0cm}
             \PARFOR{$j = 1$ to $d$}
                \State atomicAdd$\left(Y_{r_j,:}~, ~(s_j~?~A_{j, :}~:~-A_{j,:}~)\right)$
             \ENDPARFOR
         \end{algorithmic}
         \caption{CountSketch} \label{alg:countSketch}
 \end{algorithm}

 Along with its simplicity, the advantage of this CountSketch implementation is that it reduces $A \in \mathbb{R}^{d \times n}$ to $Y \in \mathbb{R}^{k \times n}$ in only $dn$ floating point reads, $d$ integer reads, $d$ boolean reads, and $dn$ floating point writes in only $dn$ floating point operations. As a consequence, the algorithm has low arithmetic intensity, and can be assumed to be memory-bound on modern architectures. Further, there are minimal storage requirements. On top of needing to store $A \in \mathbb{R}^{d \times n}$, which we will assume must always be stored, the algorithm only uses $kd$ floating point numbers, $d$ integers, and $d$ booleans. 

\section{GPU SRHT Implementation} \label{sec:srht}
As in the previous section, we begin by explicitly defining the SRHT.

\begin{definition}
    Let $\log_2d \in \mathbb{Z}$. The \emph{SRHT} is defined as the product $S = \frac{1}{\sqrt{k}}PH_dD \in \mathbb{R}^{k \times d}$, where:
    \begin{itemize}
    \item $P \in \mathbb{R}^{k \times d}$ is a row sampling matrix; i.e., each column of $P$ is a standard basis vector of dimension $k$, chosen uniformly at random,
    \item $H_d \in \mathbb{R}^{d \times d}$ is the Hadamard transform, defined recursively:
    \begin{align*}
        H_2 &= \begin{bmatrix}
            1 & 1\\
            1 & -1
        \end{bmatrix}, \quad
        H_\ell = \begin{bmatrix}
            H_{\ell/2} & H_{\ell/2}\\
            H_{\ell/2} & -H_{\ell/2}
        \end{bmatrix}, \text{ for all $\ell \in \{4, 8, 16, \dots\}$}
    \end{align*}
    \item $D \in \mathbb{R}^{d \times d}$ is diagonal, with $P(d_{i,i} = 1) = P(d_{i,i} = -1) = 1/2$.
\end{itemize}
\end{definition}

Implementing the Hadamard transform in a fast way can be accomplished using the Fast Walsh-Hadamard Transform (FWHT) \cite{FWHT}, which uses a divide-and-conquer scheme similar to the FFT. Unfortunately, the FWHT is not implemented in cuFFT, nor is any real-to-real DFT algorithm. Hence, we implemented a version using NVIDIA's CUDA samples as a guide which contain a single-vector FWHT implementation \cite{cuda-samples}, and adapted our code for multiple vectors and to maximize shared memory usage on our machine. Specifically, our implementation applies the FWHT to a matrix $A \in \mathbb{R}^{d \times n}$ column-by-column, and does so using a radix-4 FWHT, which we describe in Algorithm \ref{alg:fwhtRad4}. In the pseudocode, we neglect the details of how exactly the parallel for loop is done, where each thread is assigned to a single butterfly and the $k, b$ variables are extracted from the thread index.

\begin{algorithm}[h]
         \begin{algorithmic}[1]  
             \vspace{0cm}
             \Statex \textbf{Input/Output:}\phantom{aa} Vector $a \in \mathbb{R}^d$
             \vspace{0cm}
             \State $\texttt{stride} = d/4$
             \While{$\texttt{stride} \geq 1$}
                \State $s = \texttt{stride} * 4$
               \PARFOR{$b = 0:s:d-s,~k = 0:(\texttt{stride} - 1)$}
                        \State $i_0 = b + k + 1$, 
                        \State $i_1 = i_0 + \texttt{stride}$, $i_2 = i_0 + 2 \cdot \texttt{stride}$, $i_3 = i_0 + 3 \cdot \texttt{stride}$
                        \State $x = a_{i_0}$, $y = a_{i_1}$, $z = a_{i_2}$, $t = a_{i_3}$
                        \State $X = x + z$, $Y = y + t$, $Z = x - z$, $T = y - t$
                        \State $a_{i_0} = X + Y$, $a_{i_1} = X - Y$, $a_{i_2} = Z + T$, $a_{i_3} = Z - T$
               \ENDPARFOR
               \State $\texttt{stride}~/= 4$
             \EndWhile
         \end{algorithmic}
         \caption{Radix-4 FWHT: $\text{fwhtRad4}(a)$} \label{alg:fwhtRad4}
 \end{algorithm}

 \begin{figure}
  \centering
  \includegraphics[width=0.8\linewidth]{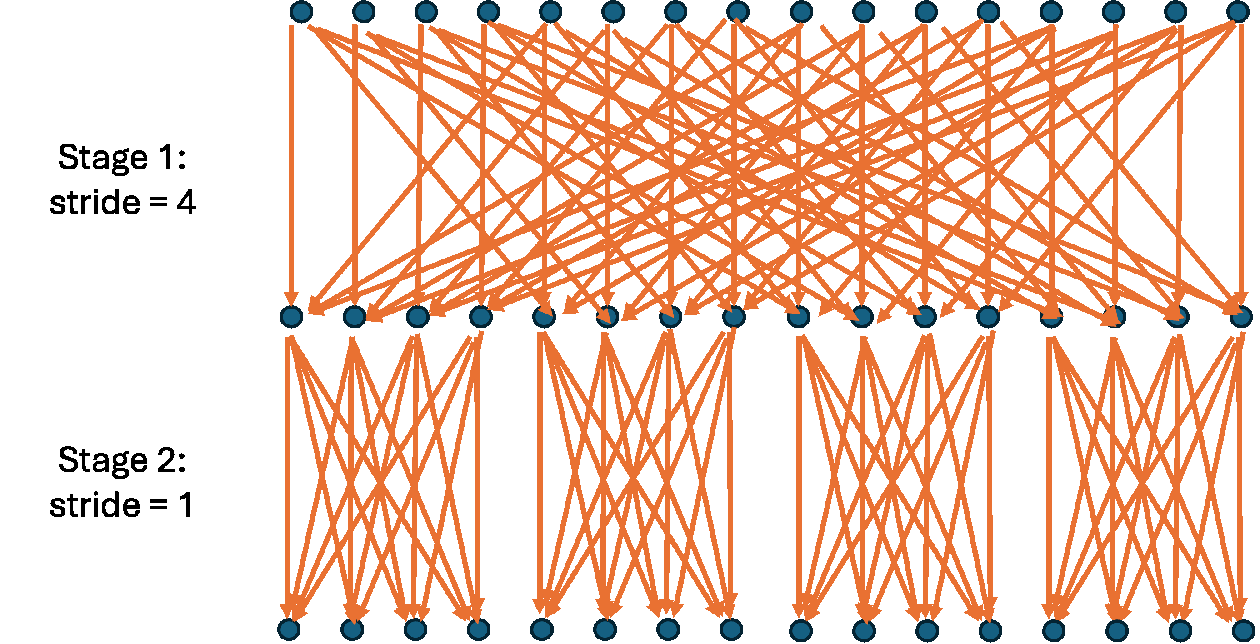}
\caption{Memory access pattern of Algorithm \ref{alg:fwhtRad4} on a vector with length $16$.}
\label{fig:butterfly}
\end{figure}

The memory access pattern used by the algorithm on a length $16$ vector is demonstrated in Figure \ref{fig:butterfly}. The figure demonstrates that as the algorithm progresses (i.e., as the stride decreases), each butterfly consists of smaller and smaller portions of the original vector. Hence, at some point these small portions eventually fit into shared memory, reducing the global memory reads necessary for the FWHT. This is precisely the way our implementation (along with the one given in the CUDA samples) works; Algorithm \ref{alg:fwhtRad4} proceeds until the stride is small enough so that the necessary portion of the input vector used in the butterfly is smaller than the maximum available shared memory on the device, at which point the device reads this portion of the input vector once into shared memory, and re-uses this for all subsequent butterfly patterns that use this portion of the vector. This can lead to substantial savings in global read/writes, since a naiive implementation requires $O(d \log d)$ read/writes on a length $d$ vector. 

Since the matrix FWHT requires $O(n d \log d)$ read/writes, we emphasized using a column-major matrix layout for the FWHT, which is more conducive to memory coalescing than row-major. Despite of the fact that memory coalescing is only possible in the later stages of the algorithm (i.e., when the stride is sufficiently small), the sheer number of read/writes make it advantageous. In contrast, an ideal parallel implementation of the row sampling and row scaling matrices $P$ and $D$, respectively, would use row-major order. However, in our experience, the cost of using column-major order for row sampling and row scaling was significantly less than the cost of the FWHT, and converting the matrix from row to column major order is significantly more expensive than just performing the row sampling and row scaling in column-major. Hence, our SRHT does all operations in column-major format.

\section{Performance Results}
In this section, we compare the performance of the CountSketch implementation described in Section \ref{sec:countsketch} on a NVIDIA H100 GPU against a naiive CountSketch implementation using a sparse matrix times dense matrix multiply (SpMM), a Gaussian sketch applied using a dense matrix matrix multiply (GeMM), and a SRHT using an implementation described in Section \ref{sec:srht}. Additionally, we incorporate our Countsketch implementation into a multisketch approach, and compare its performance against the cost of computing the Gram matrix $G = A^TA$ using GeMM. Although the Gram matrix can be computed in fewer operations using a symmetric rank-k update (SyRK), SyRK's performance is much worse in practice than GeMM.

To demonstrate the impact of the additional dimension reduction that multisketch provides over the CountSketch is useful, we used the algorithms within a sketch-and-solve least squares problem. Specifically, we compared the performance of solving a least squares problem using a sketch-and-solve approach using a Gaussian sketch, CountSketch, and the multisketching strategy (leveraging Algorithm \ref{alg:countSketch}) compared to the performance of the normal equations, which is the fastest deterministic direct least squares solver.

In some applications, the distortion that the sketch-and-solve approach introduces is unacceptable, but random sketching can also be used to quickly compute a true least squares solution. One way is to do so through an iterative method such as Blendenpik or LSRN \cite{blendenpik, lsrn}. Alternatively, one can solve it directly using a randomized Cholesky QR algorithm (also referred to in literature as randomized Householder-Cholesky QR, abbreviated as rand\_cholQR), given in Algorithm \ref{alg:randCholQR} \cite{Balabanov:2021:block, higgins2024}. We describe this algorithm in more detail in Section \ref{sec:lsqr}, and compare its performance to the aforementioned least squares methods as well. 

\subsection{Implementation Details}
All tests were done in C++ using CUDA 12.4 on a NVIDIA H100 SXM5 80GB GPU, where GeMM and SpMM operations were done leveraging the cuBLAS and cuSPARSE libraries \cite{cublas, cuda, cusparse}. The cuRAND library \cite{curand} was used to generate random numbers for the sketches. All recorded metrics were the result of an average over $100$ repeated randomly generated experiments to eliminate noise.

To achieve high performance in the CountSketch implementation given in Algorithm \ref{alg:countSketch}, memory must be read in an efficient way since the algorithm is inherently memory-bound. To coalesce memory reads, the matrix $A$ was assumed to be in row-major order so that row-wise reads of $A_{j,:}$ coincided closely with reads of $s_j$ and $r_j$ for each $j \in \{1, \dots, d\}$. Thereafter, the resulting matrix $Y$ was converted to column-major order to perform necessary tasks in cuBLAS/cuSOLVER. The one caveat to this was when multisketching was applied, the resulting $Y$ was not immediately converted to column-major order; instead, one interpreted $Y$ stored in row-major as the transpose of $Y$ stored in column-major, and computed the product $Z = GY$ where $G$ is a Gaussian matrix via $Z^T = Y^TG^T$, and subsequently transposed $Z = (Z^T)^T$. This strategy allowed us to transpose a smaller matrix, thus saving some time when converting back to column-major format. 

We used the cuSOLVER library \cite{cusolver} for our least squares experiments. The normal equations were solved by performing a Cholesky factorization on the Gram matrix $X^TX$. Sketch-and-solve algorithms were done using a QR-based solve; i.e., by using \texttt{GeQRF} followed by \texttt{OrMQR}. We avoided using \texttt{GeLS}, as this was significantly slower than using \texttt{GeQRF} and \texttt{OrMQR}.

\subsection{Sketch Performance} \label{sec:SketchPerf}
\begin{figure}
  \centering
  \includegraphics[width=\linewidth]{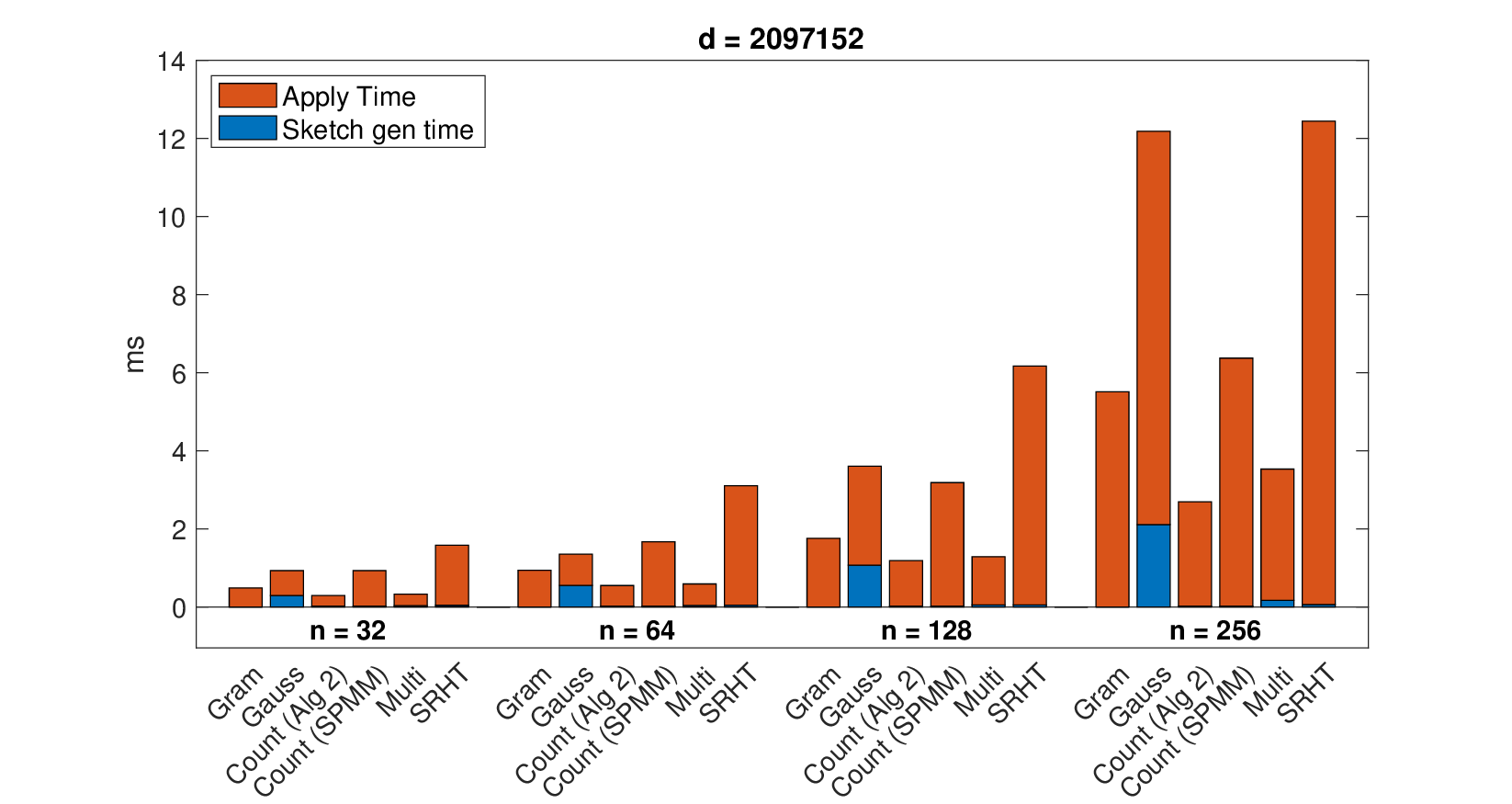}
  \label{fig:time_2e21}
  \includegraphics[width=\linewidth]{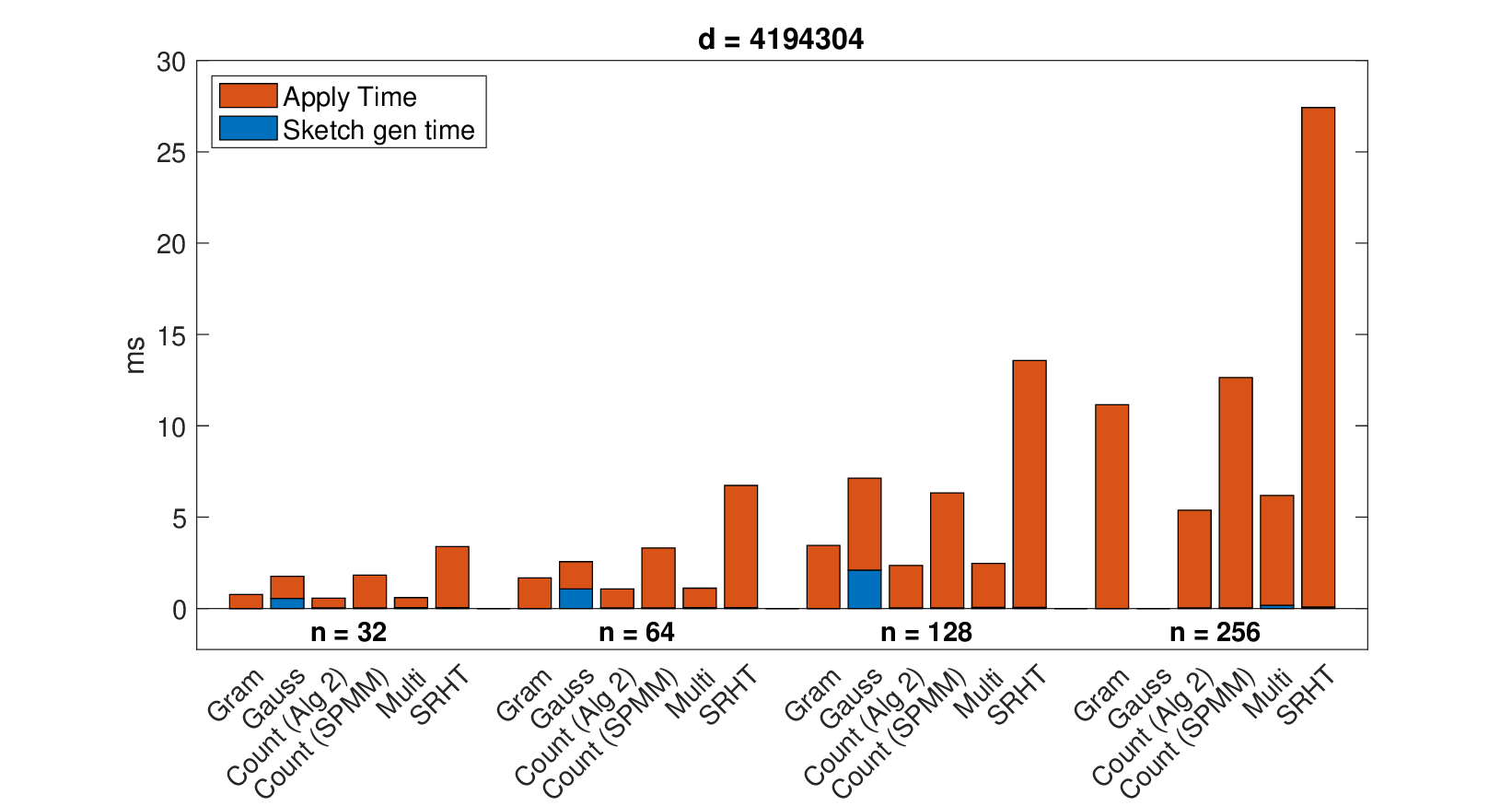}
  \label{fig:time_2e22}
  \includegraphics[width=\linewidth]{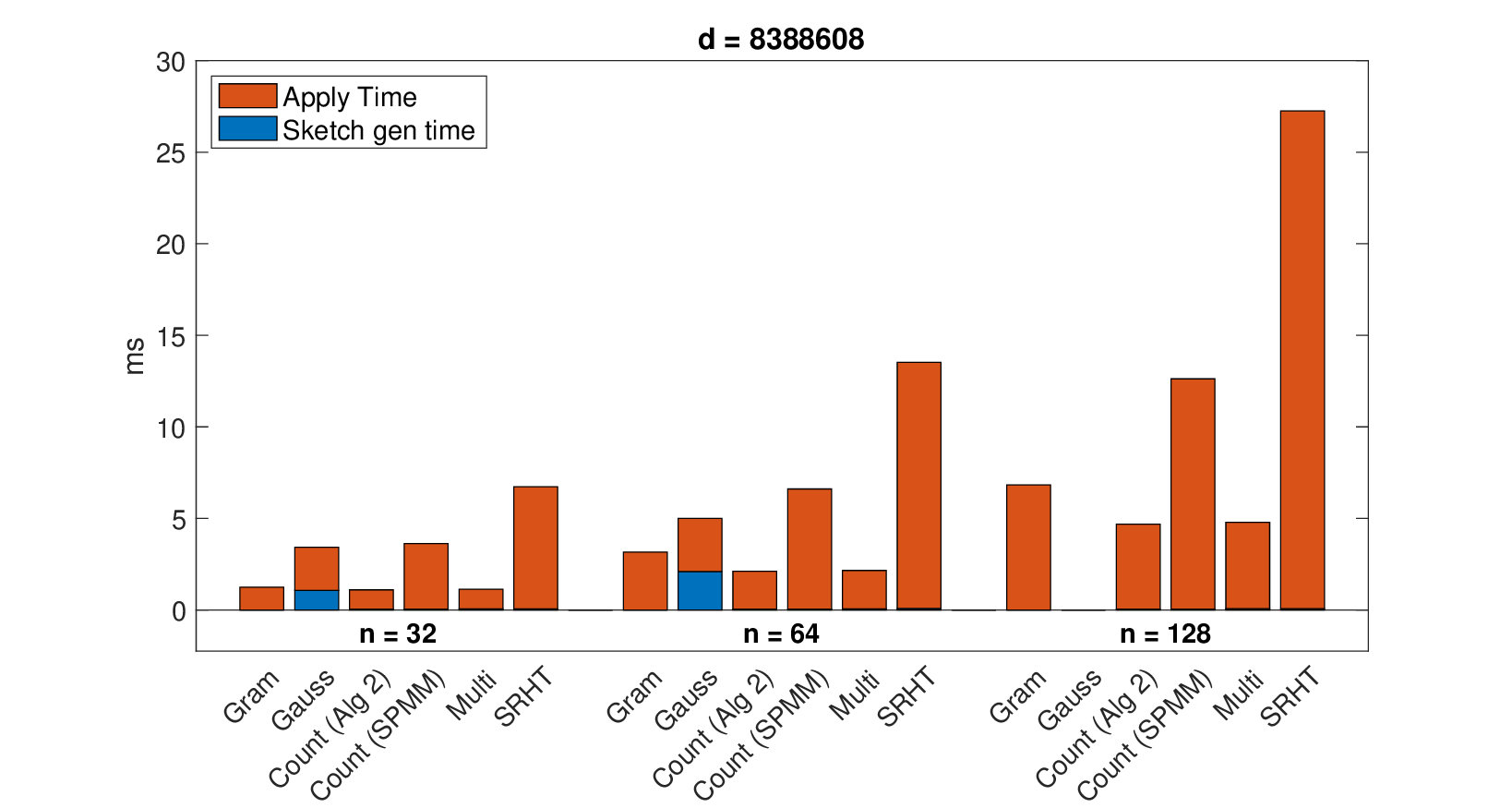}
  \label{fig:time_2e23}
%\begin{center}
%\end{center}%
\caption{Given $A \in \mathbb{R}^{d \times n}$ where $d \in \{2^{21}, 2^{22}, 2^{23}\}$ and various values of $n$, runtime breakdown of computing the Gram matrix $A^TA$ compared to applying a Gaussian sketch, CountSketch using both Algorithm \ref{alg:countSketch} and SPMM, and Multisketch using algorithm \ref{alg:countSketch}. The Gaussian bar for $d = 4194304$ and $n = 256$ is blank because the GPU ran out of memory.}
\label{fig:time_results_3}
\end{figure}

In Figure \ref{fig:time_results_3}, we fix a random matrix $A \in \mathbb{R}^{d \times n}$ with $d \in \{2^{21},2^{22},2^{23}\}$ and $n \in \{32, 64, 128, 256\}$ and report the timings required to apply a Gaussian sketch with sketch dimension $k = 2n$, a CountSketch using both the implementation given in Algorithm \ref{alg:countSketch} and an implementation using SPMM with sketch dimension $k = 2n^2$, a SRHT using the implementation described in Section \ref{sec:srht} using $k = 2n$, and a multisketch using a CountSketch (again using Algorithm \ref{alg:countSketch}) with dimension $k = 2n^2$ followed by a $2n \times 2n^2$ Gaussian matrix. 

\begin{figure}
  \centering
  \includegraphics[width=0.32\linewidth]{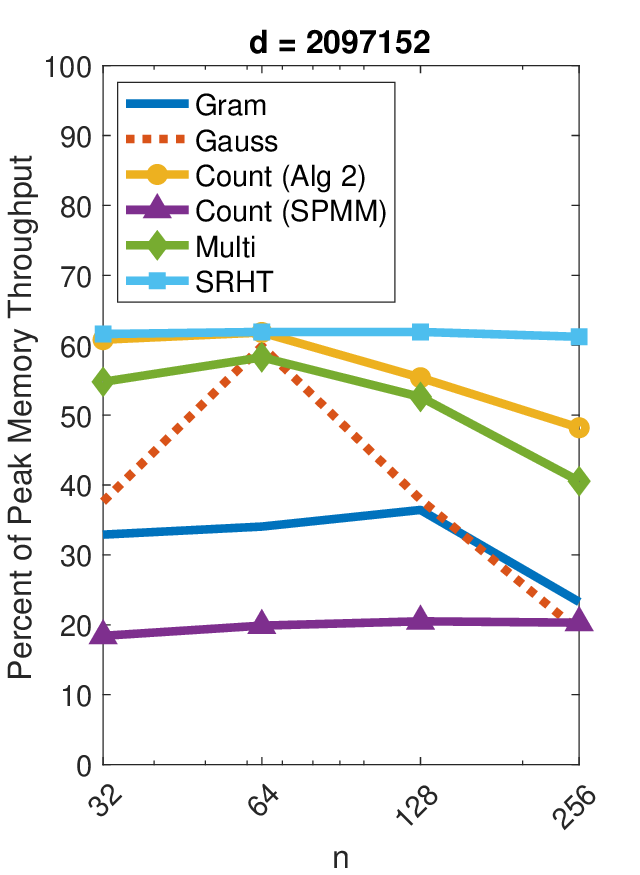}
  \label{fig:mem_2e21}
  \includegraphics[width=0.32\linewidth]{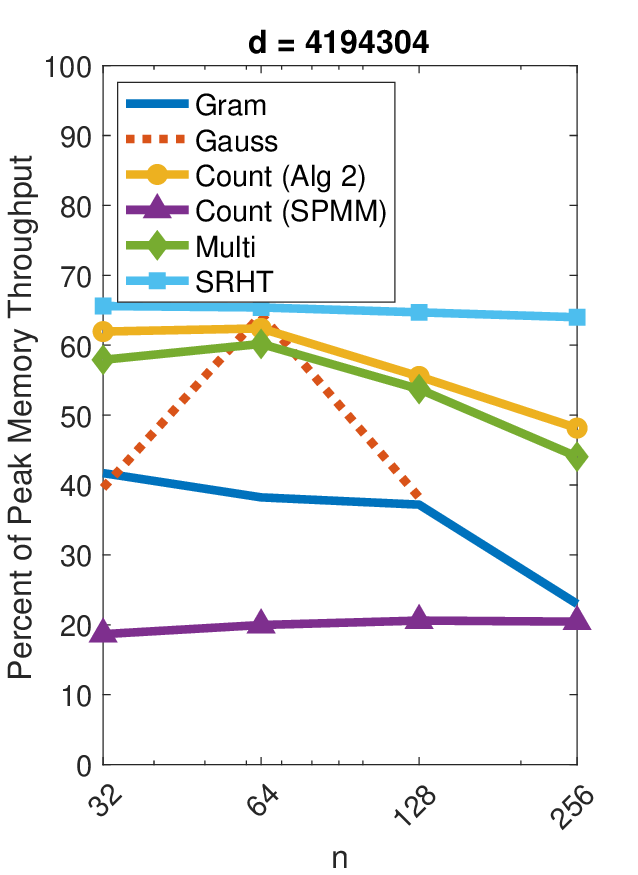}
  \label{fig:mem_2e22}
  \includegraphics[width=0.32\linewidth]{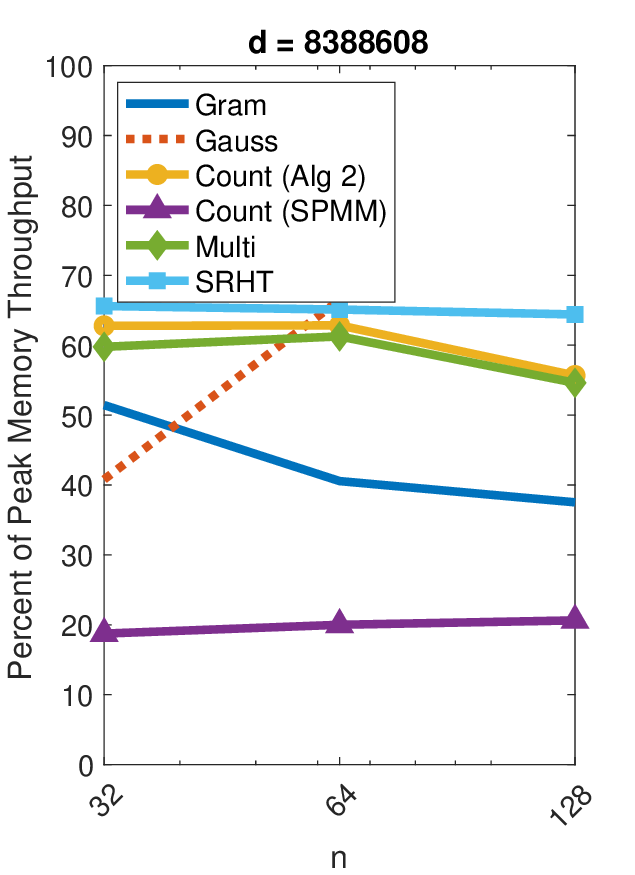}
  \label{fig:mem_2e23}
\caption{Given $A \in \mathbb{R}^{d \times n}$ where $d \in \{2^{21}, 2^{22}, 2^{23}\}$ and various values of $n$, percent of peak memory throughput of computing the Gram matrix $A^TA$ compared to applying a Gaussian sketch, CountSketch (using both algorithm \ref{alg:countSketch} and SPMM), and Multisketch using algorithm \ref{alg:countSketch}.}
\label{fig:mem_results_3}
\end{figure}

\begin{figure}
  \centering
  \includegraphics[width=0.32\linewidth]{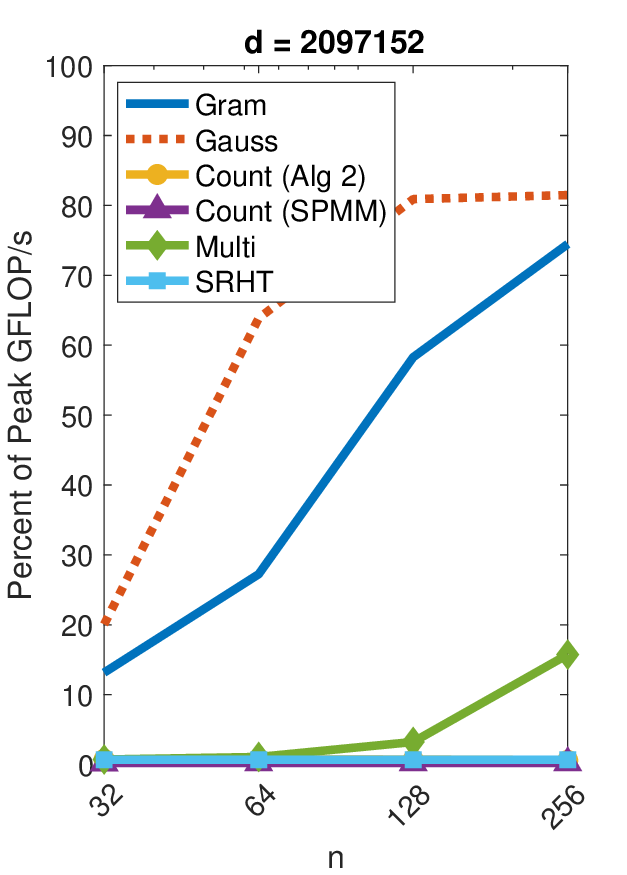}
  \label{fig:mem_2e21}
  \includegraphics[width=0.32\linewidth]{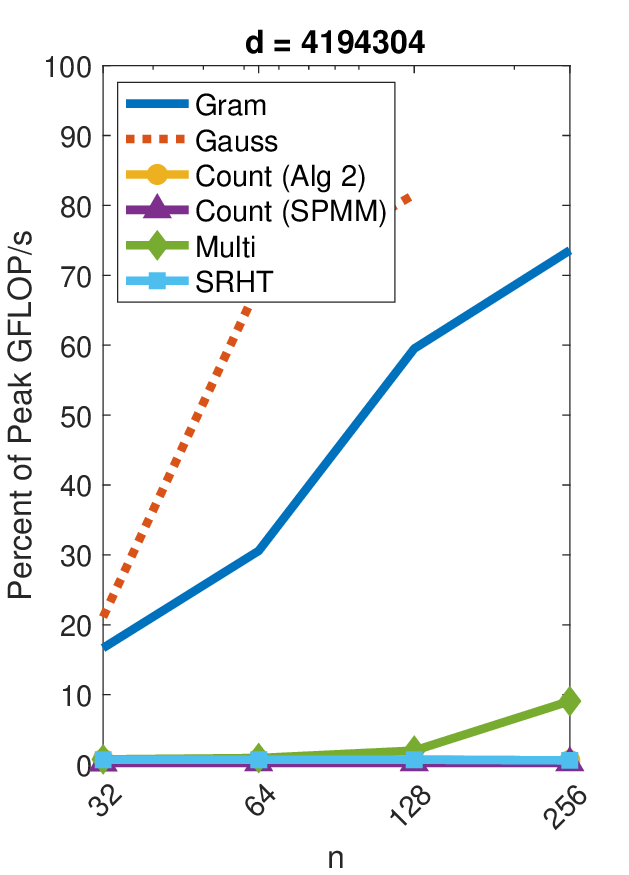}
  \label{fig:mem_2e22}
  \includegraphics[width=0.32\linewidth]{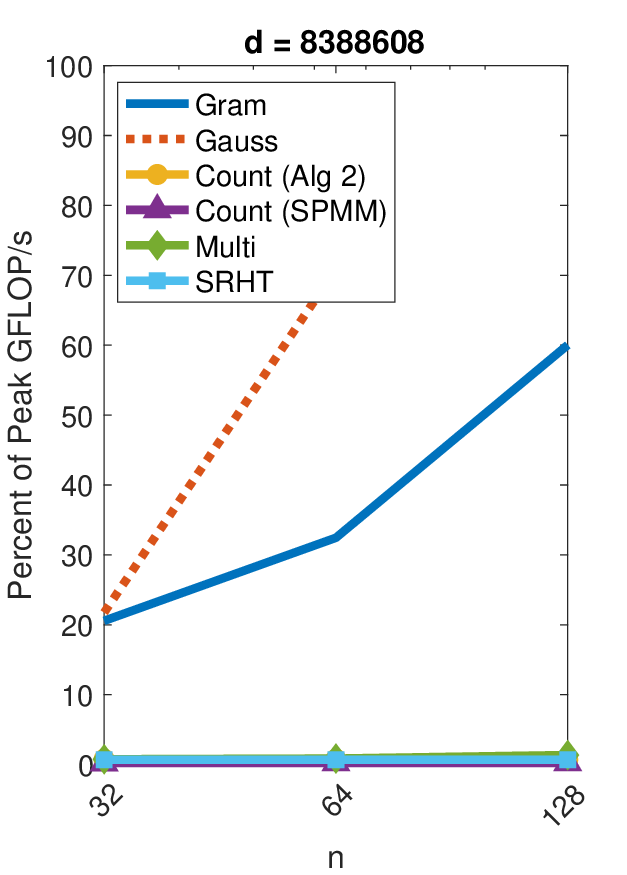}
  \label{fig:mem_2e23}
\caption{Given $A \in \mathbb{R}^{d \times n}$ where $d \in \{2^{21}, 2^{22}, 2^{23}\}$ and various values of $n$, percent of peak FLOP/s of computing the Gram matrix $A^TA$ compared to applying a Gaussian sketch, CountSketch (using both algorithm \ref{alg:countSketch} and SPMM), and Multisketch using algorithm \ref{alg:countSketch}.}
\label{fig:flop_results_3}
\end{figure}

We additionally compared these timings to the cost of computing the Gram matrix $A^TA$. This gives a natural comparison of the dominant costs of solving a least squares problem using the normal equations compared to a sketch-and-solve approach. Additionally, in Figure \ref{fig:mem_results_3}, we demonstrate the percent of peak memory throughput achieved by each of these computations. Since CountSketch is memory-bound, this figure suggests that the implementation in algorithm \ref{alg:countSketch} performs relatively well, typically hitting around $50$--$60\%$ of the memory peak, while the SPMM-based CountSketch implementation only achieves around $20\%$. There are limitations to our implementation that inherently prevent us from hitting near peak performance, such as the atomic operations, though it does stand to reason that our implementation is not quite optimal. Additionally, the SRHT implementation achieves around $60$--$70\%$ of the memory peak, which indicates that this is also implemented relatively well since the algorithm is also inherently memory-bound due to the large number of memory accesses required during the FWHT. In Figure \ref{fig:flop_results_3}, we demonstrate the percent of peak FLOP/s achieved by each computation. The algorithms that rely mostly on sparse computations to achieve dimension reduction (CountSketch, multisketch, and the SRHT) achieve poor FLOP/s, which is to be expected since these computations are memory bound.

As we can see, for sufficiently wide matrices, the CountSketch implementation provides a considerable speedup compared to computing the Gram matrix. Moreover, the multisketch technique adds minimal overhead to the CountSketch, allowing for a dimension reduction down to $k = 2n$ in significantly less time than computing the Gram matrix. Unsurprisingly, the application of a Gaussian sketch is noticably slower than computing the Gram matrix, because one performs a GeMM using a matrix that is twice as large and one has to generate $2n \times d$ i.i.d.~Gaussian random variables, which has a non-negligible overhead. In contrast, the CountSketch only needs to compute $n$ random integers and $n$ random booleans, which is significantly less expensive. Although one must also generate a dense matrix of i.i.d~Gaussians to compute the multisketch, one only needs to generate $4n^3 \ll 2nd$ Gaussians, so this overhead is also negligible. We would like to point out that there are no entries for the Gaussian sketch applied to the $2^{22} \times 256$ and $2^{23} \times 128$ problems, because storing the Gaussian caused the GPU to run out of memory for these problems.

The SRHT is generally not as competitive as the other sketching techniques. In spite of the fact that the SRHT attains better memory throughput than the CountSketch and the multisketch as indicated in Figure \ref{fig:mem_results_3}, it requires significantly more memory accesses than the other memory-bound sketch operators (CountSketch and multisketch) as indicated in Table \ref{tab:sketchProperties}. The Gaussian sketch is applied using GeMM, which is highly optimized and compute-bound which is particularly advantageous on the GPU. Further, the performance of the SRHT is limited by the performance of the FWHT, which requires significantly more synchronizations than either GeMM or the CountSketch, because of its sequential stages. 

\subsection{Least Squares Performance, Errors, and Stability}\label{sec:lsqr}

\begin{figure}
  \centering
  \includegraphics[width=\linewidth]{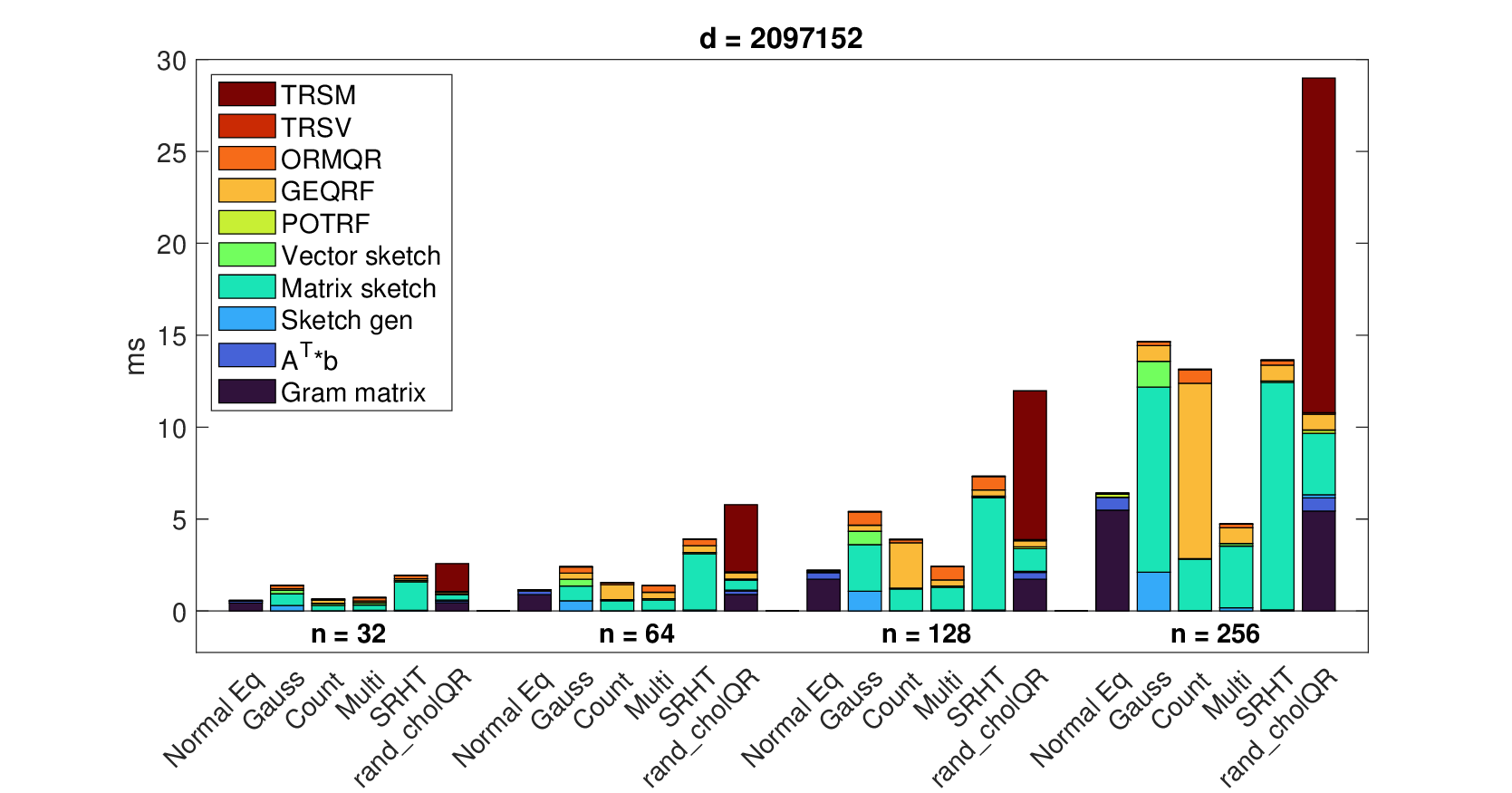}
  \label{fig:lsqr_2e21}
  \includegraphics[width=\linewidth]{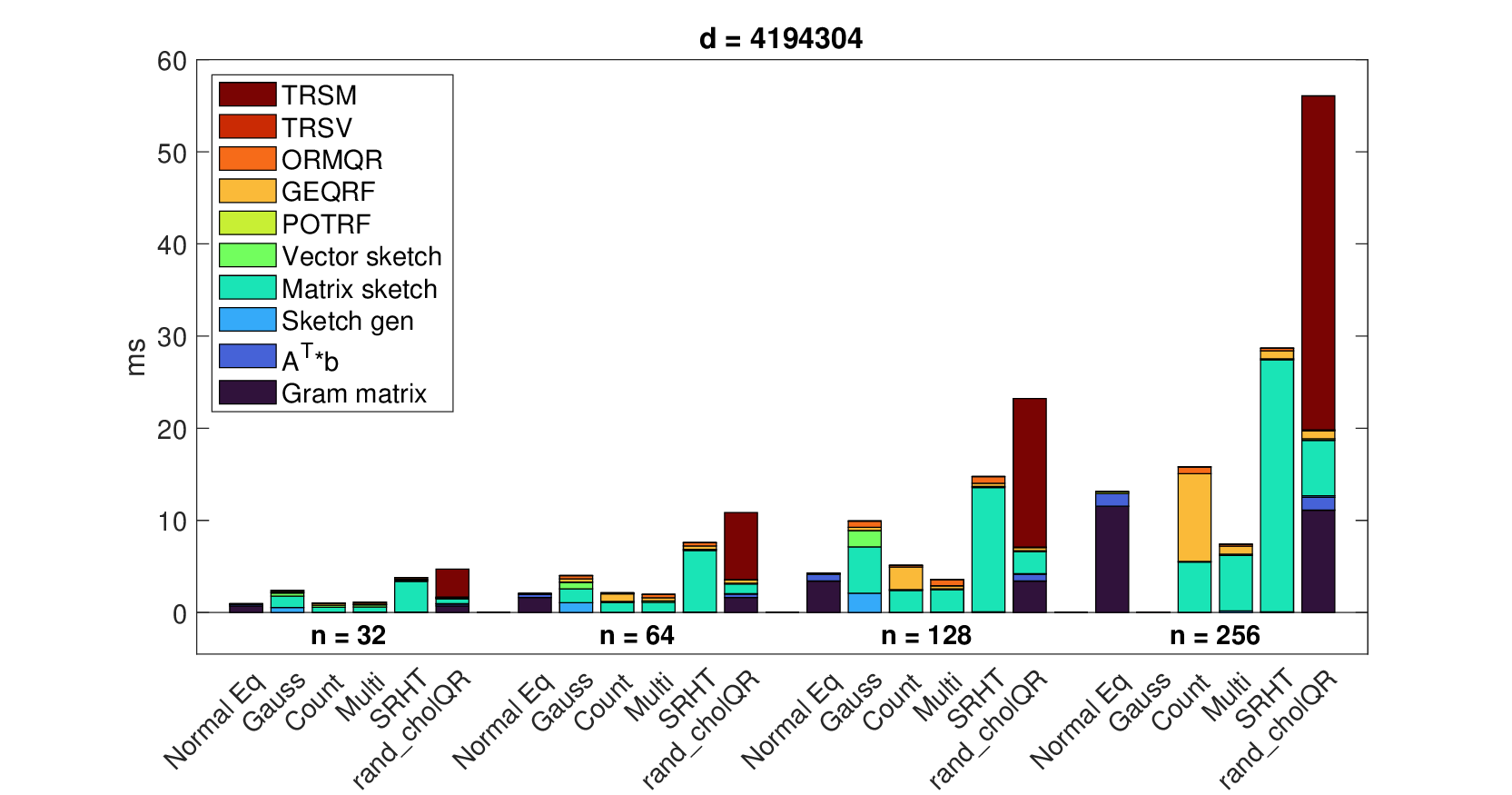}
  \label{fig:lsqr_2e22}
  \includegraphics[width=\linewidth]{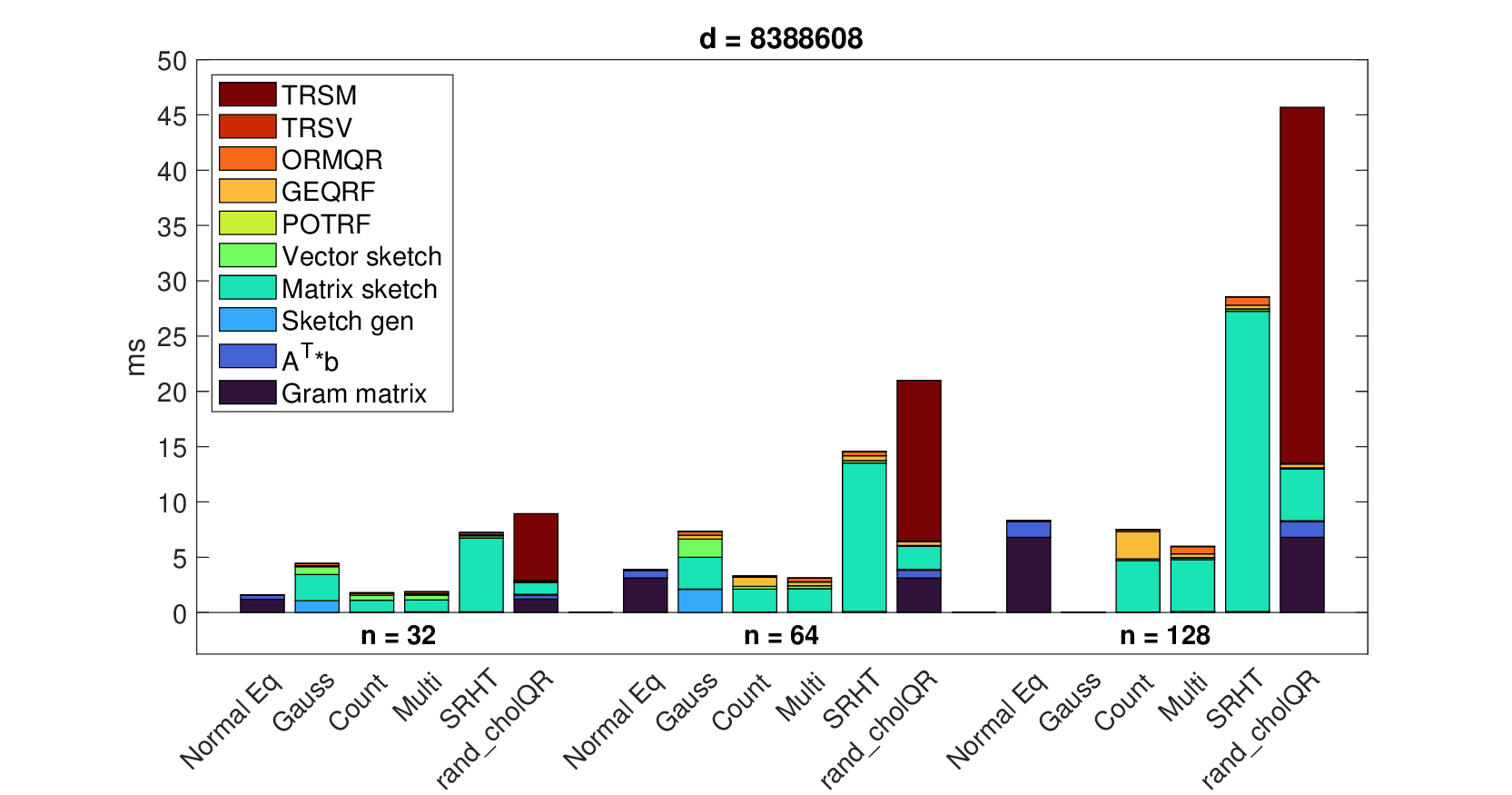}
  \label{fig:lsqr_2e23}
\caption{Given $A \in \mathbb{R}^{d \times n}$ where $d \in \{2^{21}, 2^{22}, 2^{23}\}$, $b \in \mathbb{R}^d$, and various values of $n$, runtime breakdown of solving an overdetermined least squares problem using the normal equations compared to sketch-and-solve using a Gaussian sketch, CountSketch, and Multisketch with algorithm \ref{alg:countSketch}. Blank Gaussian bars indicate the GPU ran out of memory.}
\label{fig:lsqr_time_results}
\end{figure}

In addition to testing the sketch performance, we tested the speed of using the aforementioned sketches in solving a dense least squares problem $\min \|b - Ax\|_2$, where $A \in \mathbb{R}^{d \times n}$,  $b \in \mathbb{R}^d$ for $d \in \{2^{21}, 2^{22}, 2^{23}\}$ and various values of $n$. Figure \ref{fig:lsqr_time_results} breaks down these runtimes for solving the least squares problem using the normal equations (typically the fastest direct least squares solver in practice) compared to performing a ``sketch-and-solve" algorithm with the Gaussian sketch, CountSketch, and multisketch strategy using algorithm \ref{alg:countSketch}. Again, we choose our embedding dimension for the Gaussian matrix to be $k = 2n$, and the CountSketch to be $k = 2n^2$. We did not include a standard QR-based direct solver in our performance plots, since the runtimes were substantially longer and therefore destroy the scaling of the figures.

Additionally, we compared the performance of these least squares solvers to a rand\_cholQR-based solver. Specifically, we form the true QR factorization of $A$ using rand\_cholQR, given in Algorithm \ref{alg:randCholQR}, and use this to solve the least squares problem. This algorithm was proven to be stable and accurate with high probability provided $\kappa(A) < \textbf{u}^{-1}$ where $\textbf{u}$ is unit roundoff \cite{Balabanov:2022:cholqr, higgins2024}, and was shown to perform better than a variety of other QR algorithms on a NVIDIA A100 GPU when a multisketching strategy is used \cite{higgins2024}.

\begin{algorithm}[h]
         \begin{algorithmic}[1]  
             \vspace{0cm}
             \Statex \textbf{Input:}\phantom{aa} Matrix $A \in \mathbb{R}^{d \times n}$, sketch operator $S: \mathbb{R}^d \rightarrow \mathbb{R}^k$
            \Statex {\textbf{Output:} Orthogonal $Q \in \mathbb{R}^{d\times n}$, triangular $R \in \mathbb{R}^{n \times n}$}
             \vspace{0cm}
             \State Sketch $Y = SA$
             \State Compute economy QR factorization $[\sim,R_0] = \text{qr}(Y,0)$
             \State Precondition $A$: $Q_0 = AR_0^{-1}$
             \State Compute Gram matrix: $G = Q_0^TQ_0$
             \State Compute Cholesky factorization: $R_1 = \text{chol}(G)$
             \State Compute $Q, R$: $Q = Q_0R_1^{-1}$, $R = R_1R_0$
         \end{algorithmic}
         \caption{Randomized Cholesky QR (rand\_cholQR)} \label{alg:randCholQR}
 \end{algorithm}

 In practice, one does not need the $Q$ factor to accurately solve a least squares problem using Algorithm \ref{alg:randCholQR}. Instead, one can use the approach shown in algorithm \ref{alg:randCholQRLSQR}, which performs only one TRSM operation, and is also mathematically equivalent to a recently proposed ``preconditioned normal equations" approach \cite{ipsen2025}. 
 
 \begin{algorithm}[h]
         \begin{algorithmic}[1]  
            \vspace{0cm}
             \Statex \textbf{Input:}\phantom{aa} Matrix $A \in \mathbb{R}^{d \times n}$, vector $b \in \mathbb{R}^d$, sketch operator \phantom{aaaaaaaaa}$S: \mathbb{R}^d \rightarrow \mathbb{R}^k$
            \Statex {\textbf{Output:} Solution vector $x \in \mathbb{R}^n$}
             \vspace{0cm}
             \State Sketch $Y = SA$
             \State Compute economy QR factorization $[\sim,R_0] = \text{qr}(Y,0)$
             \State Precondition $A$: $Q_0 = AR_0^{-1}$
             \State Compute Gram matrix and RHS: $G = Q_0^TQ_0$,~$z = Q_0^Tb$
             \State Compute Cholesky factorization: $R_1 = \text{chol}(G)$
             \State Compute $R = R_1R_0$
             \State Solve $y = R_1^{-T}z$
             \State Solve $x = R^{-1}y$
         \end{algorithmic}
         \caption{rand\_cholQR Least Squares} \label{alg:randCholQRLSQR}
 \end{algorithm}

It has been proven that the QR factorization in Algorithm \ref{alg:randCholQRLSQR} forms an accurate matrix decomposition of $A$ provided $\kappa(A) < \textbf{u}^{-1}$ \cite{higgins2024}, and hence the least squares solver will also be numerically stable \cite{ipsen2025}. Algorithm \ref{alg:sketchNSolve} will also be stable provided $\kappa(A) < \textbf{u}^{-1}$, because it is straightforward to show that $\kappa(SA) \leq O(1) \kappa(A)$ using the subspace embedding properties, and hence the standard Householder QR-based least squares stability analysis suggests the algorithm will be stable \cite{HighamNumAlg}. In contrast, a normal equations least squares solver will only be stable provided $\kappa(A) < \textbf{u}^{-1/2}$ \cite{TrefethenBau}. Hence, the three approaches have distinct tradeoffs: the sketch-and-solve approach is stable and can potentially solve the least squares problem the fastest, but with an $O(1)$ distortion factor introduced. The normal equations are fast and will not have any distortion, but will fail for ill-conditioned problems. The rand\_cholQR least squares is the slowest of the three, but is stable and has no distortion, and is also much faster than a standard Householder QR-based approach.

The normal equations were solved by first computing the Gram matrix $G = A^TA$, and augmenting the right hand side $y = A^Tb$ using GeMM and GeMV, respectively. The Gram matrix was then decomposed using a Cholesky factorization (POTRF) $G = R^TR$, and finally the solution was constructed with two TRSVs: $x = R^{-1}(R^{-T}y)$. In contrast, the sketch-and-solve routines were done by applying the sketch operator(s) $S$ to the coefficient matrix $W = SA$ using either GeMM in the Gaussian case or algorithm \ref{alg:countSketch} in the case of the CountSketch, and then using \texttt{GeQRF} and \texttt{OrMQR} to compute the R factor and the reflectors of the QR fatorization of $W$, and then apply the reflectors to the right hand side $b$. Finally, a TRSV was used to solve the final triangular system. Figure \ref{fig:lsqr_time_results} contains a breakdown of the runtimes of each of these experiments. To ensure stability of the normal equations, our experiments fixed $\kappa(A) = 10^2$.

Like Figure \ref{fig:time_results_3}, the bars for the Gaussian sketching on the experiments with the most columns are missing because the experiments failed because the GPU ran out of memory. For sufficiently wide matrices, the multisketch approach outperforms the normal equations solution by a substantial margin. The CountSketch is no longer the clear winner, because it takes a large performance hit during the \texttt{GeQRF} phase, since it must perform \texttt{GeQRF} on a larger problem than the implementations using the multisketch and the Gaussian sketch. In spite of this, the CountSketch-based sketch-and-solve algorithm sometimes still outperforms the normal equations.

\begin{figure}
  \centering
  \includegraphics[width=0.32\linewidth]{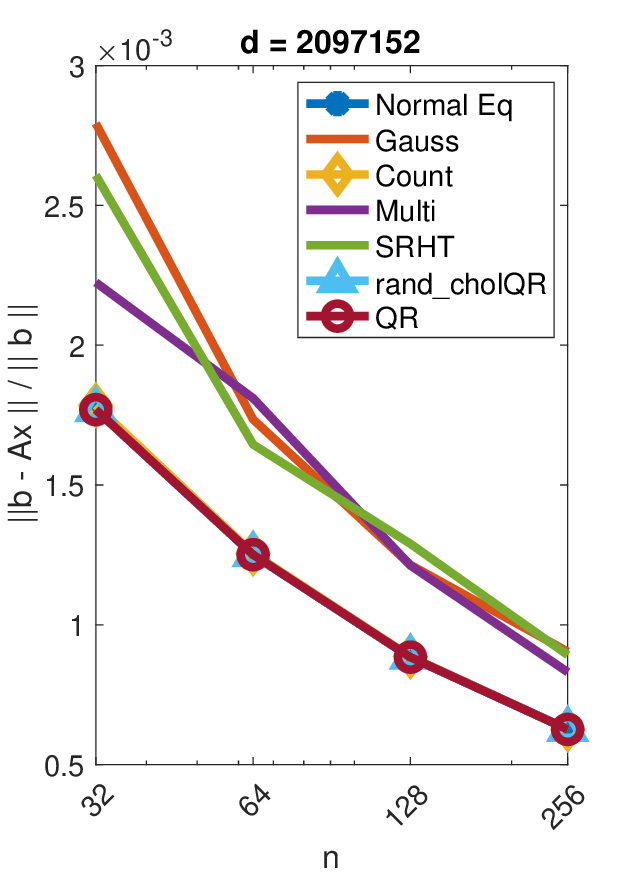}
  \includegraphics[width=0.32\linewidth]{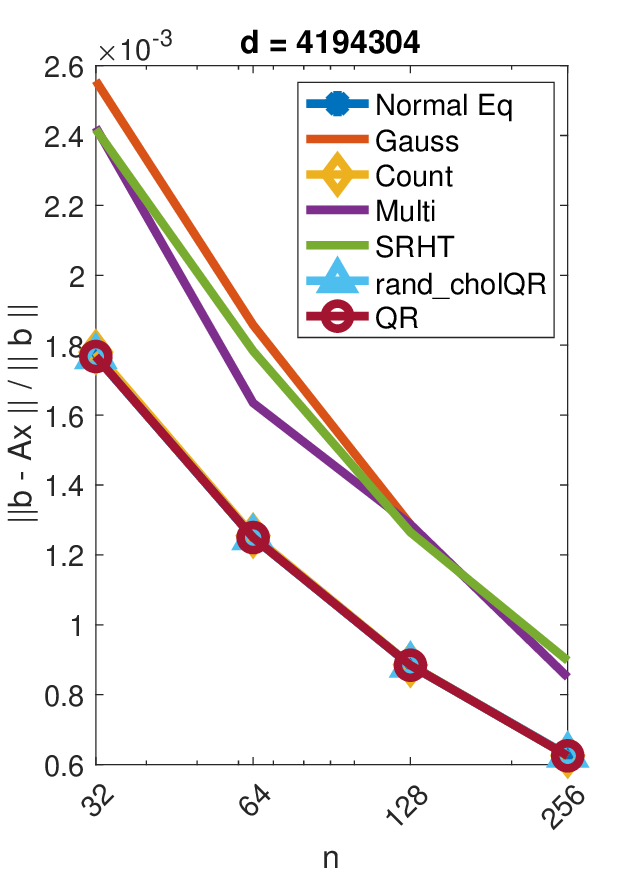}
  \includegraphics[width=0.32\linewidth]{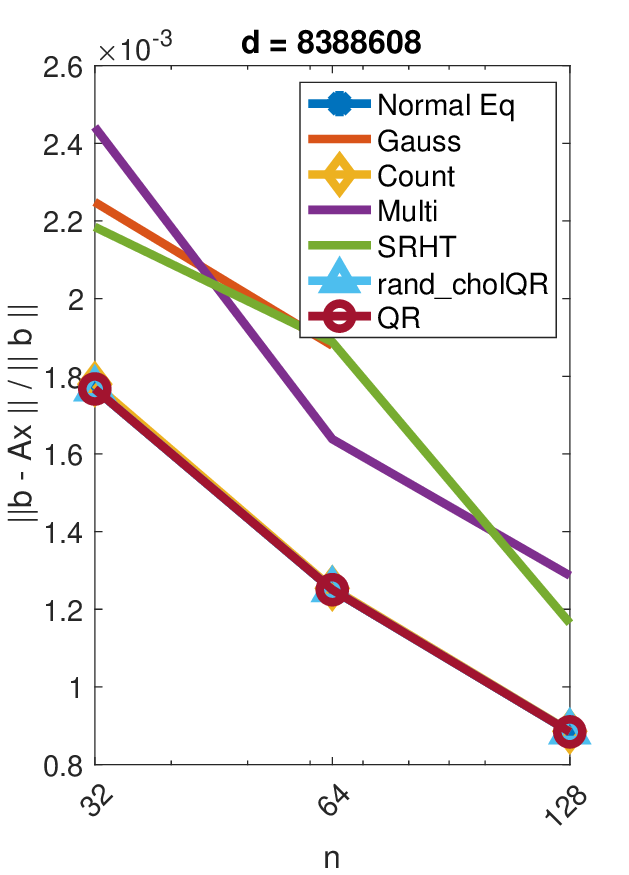}
\caption{Given $A \in \mathbb{R}^{d \times n}$ where $d \in \{2^{21}, 2^{22}, 2^{23}\}$ and various values of $n$, relative least squares residuals $\|b - Ax\|_2 / \|b\|_2$ of the ``easy" problem,}
\label{fig:lsqr_err1}
\end{figure}

\begin{figure}
  \centering
  \includegraphics[width=0.32\linewidth]{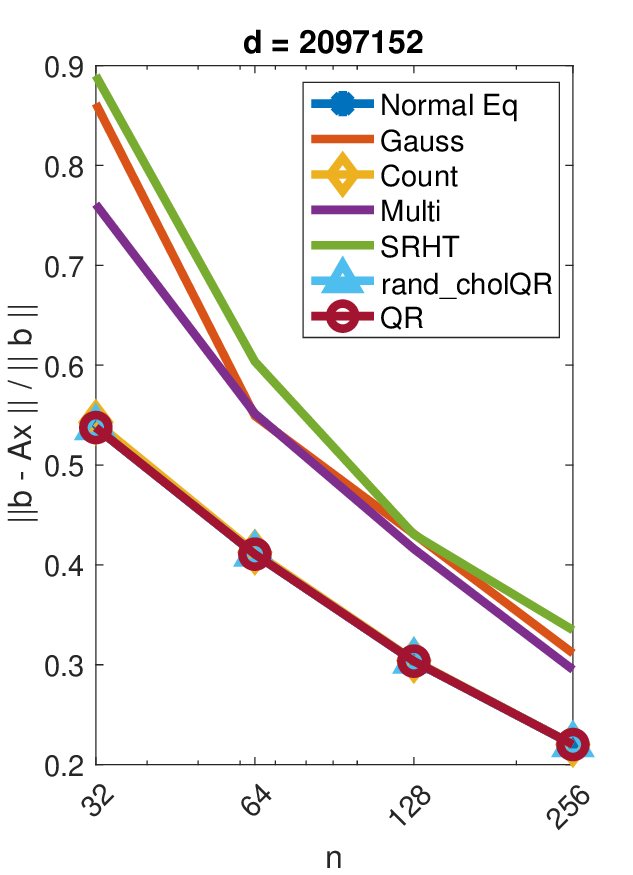}
  \includegraphics[width=0.32\linewidth]{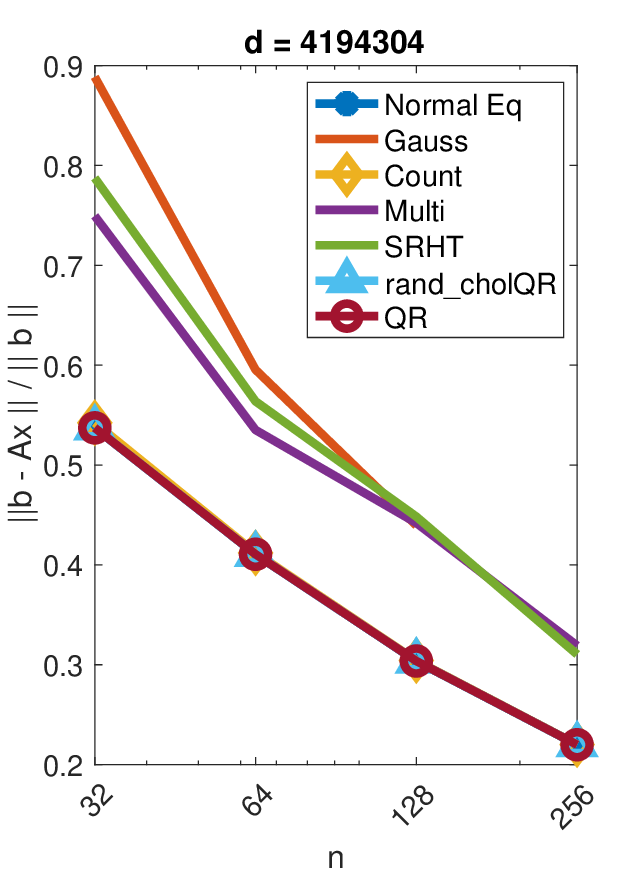}
  \includegraphics[width=0.32\linewidth]{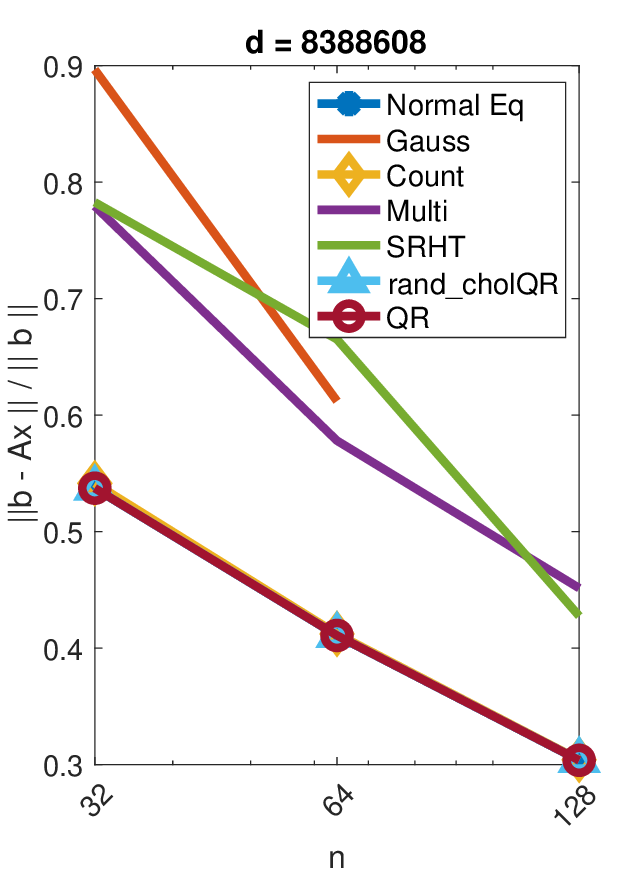}
\caption{Given $A \in \mathbb{R}^{d \times n}$ where $d \in \{2^{21}, 2^{22}, 2^{23}\}$ and various values of $n$, relative least squares residuals $\|b - Ax\|_2 / \|b\|_2$ of the ``hard" problem,}
\label{fig:lsqr_err2}
\end{figure}

Although the multisketch algorithm is faster than the normal equations for most of our experiments, its solution introduces a non-negligible distortion in the least squares relative residual $\| b - Ax \|_2 / \| b \|_2$ due to the subspace embedding properties discussed in Section \ref{sec:intro}. In Figures \ref{fig:lsqr_err1}--\ref{fig:lsqr_err2}, we demonstrate the relative residual for each experiment, first on a least squares problem designed to have a low relative residual, and another on a problem with much higher relative residual. To do so, we set $e = [1, \dots, 1]^T \in \mathbb{R}^n$ and set the right hand side $b = Ae + \eta$, where $\eta \in \mathbb{R}^d$ has $\eta_i \in \mathcal{N}(\mu, \sigma^2)$, for $\mu \in \{0, 3\}$ and $\sigma^2 \in \{0.01,2\}$ to simulate an ``easy" problem with low noise, and a ``hard" problem with higher noise.

The least squares residuals for the ``sketch-and-solve" approach generally are proportional to the true residual, which is measured here by the normal equations, which is valid since the problem is very well-conditioned. However, as proven in \cite{higgins2024}, the QR factorization formed by the ``sketch-and-solve" approach is stable provided $\kappa(A) \lessapprox \textbf{u}^{-1}$ where $\textbf{u}$ is unit roundoff. In contrast, the normal equations only maintain stability for $\kappa(A) \lessapprox \textbf{u}^{-1/2}$, as the Gram matrix begins to lose numerical rank beyond that point. 

\begin{figure}
  \centering
  \includegraphics[width=0.6\linewidth]{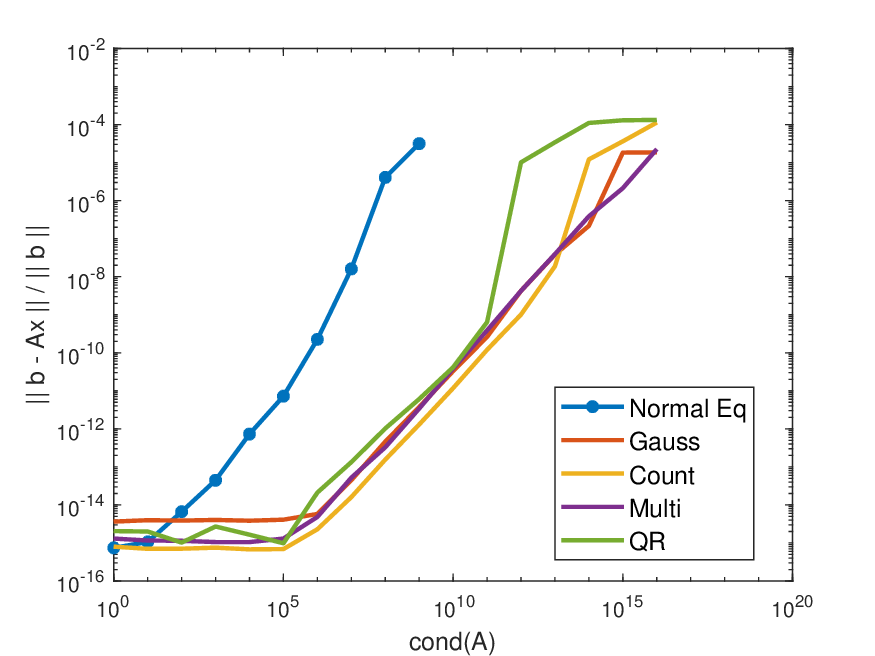}
\caption{Given $A \in \mathbb{R}^{d \times n}$ where $d = 2^{17}$ and $n = 16$, we solve $\min \| b - Ax \|_2$ using $b = Ae$, where $e$ is a vector of all ones, and $A$ has a varying condition number.}
\label{fig:sensitivity}
\end{figure}

To demonstrate the impact of this, we show the sensitivity of the least squares residual $|| b - Ax||_2 / || b ||$ to the condition number of the input matrix $A$ using the normal equations, the sketch-and-solve algorithms, and a QR-based least squares solve in Figure \ref{fig:sensitivity}. The specific problem was chosen so that an exact solution exists in exact arithmetic; specifically, we chose the right hand side $b = Ae$.

The behavior of the sketch-and-solve  approach is similar to that of the QR-based solver, while the normal equations are significantly less stable, and fail for $\kappa(A) > 10^8$. Thus, the multisketched sketch-and-solve least squares solver offers the opportunity to solve a least squares problem faster than the normal equations in a much more stable manner, at the cost of a constant $O(1)$ distortion factor.

\section{Considerations for Distributed Environments}

All of the experiments in this paper focused on the sketch performance on a single GPU. However, many of the conclusions translate to the distributed case as well. Suppose our matrix $A$ is distributed across $p$ processors in a block-row format, where $A^{(i)}$ is the block of $A$ stored on processor $i$. If we wish to apply a Gaussian sketch $G$ to $A$, $G$ consists of i.i.d.~Gaussian random variables, we can simply generate a Gaussian matrix $G^{(i)}$ on the $i^{th}$ process, apply it to $A^{(i)}$ and reduce to a single process at the end. In other words, we define $G = [G^{(1)}, G^{(2)}, \dots, G^{(p)}]$, and then $GA = \sum_{i = 1}^p G^{(i)}A^{(i)}$.

Similarly, by Definition \ref{def:countSketch}, one can interpret the CountSketch as a sparse matrix with exactly one $\pm 1$ per column. Hence, a similar strategy can be used to apply an entire CountSketch $C$ to $A$ by generating a smaller CountSketch $C^{(i)}$ on each process and again reducing at the end. In other words, define $C = [C^{(1)}, C^{(2)}, \dots, C^{(p)}]$, and then $CA = \sum_{i = 1}^p C^{(i)}A^{(i)}$. In a distributed setting, the CountSketch will have a higher communication cost than the Gaussian, since each $C^{(i)}A^{(i)}$ term is larger than the $G^{(i)}A^{(i)}$ terms, due to the larger embedding dimension of the CountSketch. However, the performance of the CountSketch per-process is significantly better than that of the Gaussian, as indicated in Section \ref{sec:SketchPerf}.

Since the Gaussian matrix $G_{ms} \in \mathbb{R}^{2n \times 2n^2}$ used in the multisketch is small, one can cast it to each process, and then apply a smaller multisketch on each process, and reduce. Specifically,
\begin{align*}
    G_{ms}CA &= G_{ms}\sum_{i = 1}^p C^{(i)}A^{(i)}
    = \sum_{i = 1}^pG_{ms} C^{(i)}A^{(i)}.
\end{align*}
Conveniently, the communication cost of the multisketch will be the same as that of the Gaussian, due to the fact that they have the same embedding dimension. Because it shares the same communication cost and has significantly better performance per-process (see Section \ref{sec:SketchPerf}), the multisketch will almost certainly outperform the Gaussian in a distributed setting as well.

Due to the memory access pattern of the FWHT used in the SRHT, distributed implementations of the standard SRHT are significantly more challenging and likely lead to performance degredation due to increasing communication costs. However, one can instead use a block SRHT \cite{balabanovBlockSRHT}, which is defined by assigning a SRHT to blocks of columns of the sketch. This approach has been proven to be a subspace embedding provided $k = O(n \log n)$, like the original SRHT, but is significantly more efficient in distributed environments. However, since the performance per-process of the SRHT is much worse than that of the CountSketch and the multisketch, and it incurs a higher communication cost than the multisketch due to its larger embedding dimension, one can also expect worse performance of this sketch compared to the multisketch.

\section{Conclusions and Future Work}
In this work, we used a simple but efficient CountSketch GPU implementation in Algorithm \ref{alg:countSketch} that allows one to apply either a CountSketch or multisketch to $A \in \mathbb{R}^{d \times n}$ significantly faster than one can compute the Gram matrix or a Gaussian sketch of the same matrix. Further, we demonstrate that using this implementation can solve a least squares problem up to $77\%$ faster than the normal equations (in the case of $d = 2^{22}$ and $n = 256$) with much better stability, at the cost of maginifying the relative residual norm $\| b - Ax \|_2 / \| b \|_2$ by an $O(1)$ factor.

Future work includes modifying Algorithm \ref{alg:countSketch} to build the CountSketch matrix on the fly using a hash-based strategy, as was intended in the original CountSketch paper \cite{CountSketchOriginal}. This presumably would increase the sketch application + generation time, but would make the method more amenable to streaming applications. We are also interested in testing other multisketching implementations outside of simply using a CountSketch with a Gaussian sketch, such as using a CountSketch with a SRHT.

\begin{acks}
    Sandia National Laboratories is a multimission laboratory managed and operated by National Technology and Engineering Solutions of Sandia, LLC., a wholly owned subsidiary of Honeywell International, Inc., for the U.S. Department of Energy’s National Nuclear Security Administration under contract DE-NA-0003525.
\end{acks}

\bibliographystyle{ACM-Reference-Format}
\bibliography{sample-base}

\end{document}